\documentclass{amsart}
\usepackage{amsmath}
\usepackage[latin1]{inputenc}
\usepackage{amssymb, amsmath}

\theoremstyle{plain}
\newtheorem{theorem}{Theorem}[section]
\newtheorem{proposition}[theorem]{Proposition}
\newtheorem{fact}[theorem]{Fact}
\newtheorem{lemma}[theorem]{Lemma}
\newtheorem{corollary}[theorem]{Corollary}
\newtheorem*{claim}{Claim}

\theoremstyle{definition}
\newtheorem{definition}[theorem]{Definition}
\newtheorem{remark}[theorem]{Remark}
\newtheorem*{question}{Question}

\def\claim{\begin{claim}}
\def\eclaim{\end{claim}}
\def\defi{\begin{definition}}
\def\edefi{\end{definition}}
\def\teo{\begin{theorem}}
\def\eteo{\end{theorem}}
\def\fet{\begin{fact}}
\def\efet{\end{fact}}
\def\cor{\begin{corollary}}
\def\ecor{\end{corollary}}
\def\lem{\begin{lemma}}
\def\elem{\end{lemma}}
\def\preg{\begin{question}}
\def\epreg{\end{question}}
\def\rem{\begin{remark}}
\def\erem{\end{remark}}
\def\prova{\begin{proof}}
\def\eprova{\end{proof}}

\def\prop{\begin{proposition}}
\def\eprop{\end{proposition}}

\def\Ind{\setbox0=\hbox{$x$}\kern\wd0\hbox to 0pt{\hss$\mid$\hss}
\lower.9\ht0\hbox to 0pt{\hss$\smile$\hss}\kern\wd0}
\def\Notind{\setbox0=\hbox{$x$}\kern\wd0\hbox to 0pt{\mathchardef
\nn=12854\hss$\nn$\kern1.4\wd0\hss}\hbox to
0pt{\hss$\mid$\hss}\lower.9\ht0 \hbox to
0pt{\hss$\smile$\hss}\kern\wd0}
\def\ind{\mathop{\mathpalette\Ind{}}}
\def\nind{\mathop{\mathpalette\Notind{}}}

\newcommand{\sub}{\subseteq}
\def\bdd{\mathrm{bdd}}
\def\dcl{\mathrm{dcl}}
\def\acl{\mathrm{acl}}
\def\cb{\mathrm{Cb}}
\def\tp{\mathrm{tp}}

\def\Aut{\mathrm{Aut}}
\def\pwt{\mathrm{pwt}}
\def\w{\mathrm{w}}

\def\D{\mathrm{D}}

\title{On omega-categorical simple theories}
\date{October 28, 2011}
\author{Daniel Palac\'\i n}\thanks{The author was partially supported by research project MTM 2008-01545 of the Spanish government and research project 2009SGR 00187 of the Catalan government. }%\newline This work was partially done while during January 2010 the author was visiting the Logic group of the University of East Anglia in Norwich. The author would like to thank Enrique Casanovas and David M. Evans for fruitful discussions.}
\keywords{simple; CM-trivial; $\omega$-categorical; low; weight; strong stable forking}
\subjclass[2000]{03C45}

\address{Universitat de Barcelona; Departament de L\`ogica, Hist\`oria
i Filosofia de la Ci\`encia, Montalegre 6, 08001 Barcelona, Spain}
\email{dpalacin@ub.edu}

\begin{document}

\begin{abstract}  In the present paper we shall prove that countable $\omega$-categorical simple CM-trivial theories and countable $\omega$-categorical simple theories with strong stable forking are low. In addition, we observe that simple theories of bounded finite weight are low. \end{abstract}

\maketitle\section{Introduction}

Buechler \cite{bue} and Shami \cite{shami} introduced the class of simple {\em low} theories which includes all currently known natural examples of simples theories. Namely, stable theories and supersimple theories of finite $\D$-rank are low. Moreover, as we will observe, simple theories of bounded finite weight are also low. However, Casanovas and Kim \cite{cas-kim} proved the existence of a supersimple nonlow theory. Roughly speaking, lowness implies that dividing for a formula is type-definable; as a consequence Buechler and Shami, independently, solved one of the most important conjectures in simple theories for simple low theories, i.e., Lascar strong types and strong types coincide for such theories. Moreover, Ben-Yaacov, Pillay, and Vassiliev \cite{ben-pil-vas}, in order to generalize Poizat's {\em Belles Paires} to their so called {\em Lovely Pairs}, observe that lowness is not far from the right simple analogue of the non finite cover property. Therefore, it turns out that lowness seems to be a very natural notion for simple theories.
Casanovas and Wagner \cite{cas-wag} introduced another interesting class of simple theories which contains all supersimple theories and simple low theories; such theories are called {\em short}. Indeed, this subclass was already considered in \cite{cas1} where Casanovas obtains a simple nonshort theory. It is worth remarkable that Casanovas' example and Casanovas and Kim's example are both one-based, i.e., every canonical base $\cb(a/A)$ is contained in $\bdd(a)$.

As the title of the paper suggests we will be concerned with the $\omega$-categorical simple framework, where much less is known. In \cite{cas-wag}, Casanovas and Wagner showed that shortness and lowness coincide in this
%framework
context and they asked the following question:
\preg Is every $\omega$-categorical simple theory low? \epreg

It is well-known that $\omega$-categorical simple one-based theories are supersimple and so, they are low. The present paper is devoted to answer affirmatively this question under the assumption of CM-triviality, a geometric property introduced by Hrushovski \cite{hrush} which generalizes one-basedness. At the time of writing, all known examples of $\omega$-categorical simple theories are CM-trivial. In particular, those obtained from a Hrushovski construction with a standard predimension are CM-trivial \cite{evans,Wagner-book}, and it seems a significant variation on the construction would be required in order to produce $\omega$-categorical simple theories which are not CM-trivial. Moreover, as we will see, in order to obtain an $\omega$-categorical simple nonlow theory the non-forking independence cannot come from finite sets.

Another approach to our question is via stability of forking. Kim and Pillay introduced a strong version of the stable forking conjecture \cite{kim-pil}. Even though they show that any completion of the theory of pseudofinite fields does not satisfy this strong version of stable forking, we think interesting to study such property in the $\omega$-categorical framework. We will show that $\omega$-categorical simple theories with strong stable forking are low. In fact, Kim and Pillay proved that one-based theories with elimination of hyperimaginaries have strong stable forking, and so do the mentioned examples of simple nonlow theories due to Casanovas and Casanovas and Kim.\vspace{0.3cm}

Part of the work was done in January 2010 during a research stay in the University of East Anglia, Norwich. I thank Enrique Casanovas and David Evans for helpful discussions.

\section{Preliminaries}

We will consider a complete first-order theory $T$ (with infinite models) in a language $L$ whose monster model is denoted by $\mathfrak{C}$. As usual, tuples and sets of
parameters will live in $\mathfrak{C}^{eq}$, and given any two tuples $a,b$ and any set of parameters $A$, we shall write $a\equiv_A b$ whenever $a$ and $b$ have the same type over $A$. We assume the reader is familiarized with the general theory of simplicity and hyperimaginaries; otherwise we recommend \cite{cas2, Wagner-book}.

%In this section we present some definitions and results regarding subclasses of simple theories.

\subsection{Lowness} Buechler \cite{bue} and Shami \cite{shami} introduced lowness using $\D$-ranks. Following \cite{ben-pil-vas}, a formula $\varphi(x;y)\in L$ is {\em low} if there is some $k<\omega$ such that for every indiscernible sequence $(a_i:i<\omega)$ the following holds: if $\{\varphi(x,a_i):i<\omega\}$ is inconsistent, then it is
$k$-inconsistent. However, for our purposes it is better to deal with dividing chains since it is easier (at least for us) to understand the relation between lowness, shortness and simplicity.

\defi Let $\alpha$ be an ordinal. A formula $\varphi(x;y)\in L$
divides $\alpha$ times if there is a sequence $(a_i:i<\alpha)$ in the
monster model such that $\{\varphi(x,a_i):i<\alpha\}$ is consistent
and $\varphi(x,a_i)$ divides over $\{a_j:j<i\}$ for all
$i<\alpha$. A such sequence $(a_i:i<\alpha)$ is called a
dividing chain of length $\alpha$. \edefi

\fet\cite[Remark 2.2]{cas1} A theory is simple iff no formula divides $\omega_1$
times iff no formula divides $\omega$ times with respect to some
fixed $k<\omega$.\efet
%A formula $\varphi(x;y)\in L$ is said to divide arbitrarily often if it divides $n$ times for all $n<\omega$.

\defi A formula $\varphi(x;y)\in L$ is {\em short} if does not divide infinitely many times; and it is {\em low} if there is some $n<\omega$ such that it does not divide $n$ times. We say a theory is short (low) if all formulas are short (low).\edefi

In \cite{cas-wag} it is remarked that our definition coincide with \cite{bue,shami}. In \cite{ben-pil-vas} it is showed that all these notions of lowness are the same for a simple theory, and in particular it is remarked that every stable formula is low \cite[Remark 2.2]{ben-pil-vas}. In addition, in \cite{ben-pil-vas} it is proved the equivalence between lowness and the $\emptyset$-type-definability of dividing over arbitrary sets of parameters.

\fet\cite[Lemma 2.3]{ben-pil-vas}\label{fact-type-def} Assume $T$ is simple. Then, $\varphi(x,y)\in L$ is low iff the relation on $y,z$ `$\varphi(x,y)$ divides over $z$' is $\emptyset$-type-definable ($z$ may be of infinite length). \efet

It is clear from our definition that low theories and supersimple theories are short, and short theories are simple.  In the $\omega$-categorical simple context we know the following.

\fet\cite[Proposition 19]{cas-wag}\label{Fact-indiscernibles} In an $\omega$-categorical theory a short formula is low. Moreover, if a formula is nonshort, then there is an indiscernible sequence witnessing this. \efet

In fact, the proof given by Casanovas and Wagner shows more:

\lem\label{indisc-seq} Let $T$ be $\omega$-categorical and let
$\varphi(x,y)$ be a nonlow formula. Then, there is a tuple $c$ (of the right length)
and some $c$-indiscernible sequence $(a_i:i<\omega)$ such that
$c\models\bigwedge_{i<\omega}\varphi(x,a_i)$ and also
$\varphi(x,a_i)$ divides over $\{a_j:j<i\}$ for all $i<\omega$.\elem
\prova We offer a proof for convenience. Since $T$ is $\omega$-categorical and $\varphi(x,y)$ is
nonlow, $\varphi(x,y)$ is nonshort by Fact \ref{Fact-indiscernibles}. Let
$(a_i:i<\omega)$ be a sequence exemplifying that $\varphi(x,y)$ divides $\omega$
times. In particular, there is some
$c\models\bigwedge_{i<\omega}\varphi(x,a_i)$. By
$\omega$-categoricity we may assume that $a_0\equiv_ca_i$ for all
$i<\omega$. Then, by $\omega$-categoricity and Ramsey's Theorem,
there is an infinite $2$-indiscernible over $c$ subsequence.
Iterating this process we infer that for all $n\geq 1$ there is an
infinite $n$-indiscernible over $c$ subsequence. By compactness, the
limit type  $q$ of these subsequences exists and if
$(b_i:i<\omega)\models q$, then it is an indiscernible sequence over
$c$. Moreover, for every $k<\omega$ there is a sequence $(n_i:i\leq
k)$ such that $(b_i:i\leq k)\equiv_c(a_{n_i}:i\leq k)$. Thus, since
$\varphi(x,a_{n_k})$ divides over $\{a_{n_i}:i<k\}$, so does
$\varphi(x,b_k)$ over $\{b_i:i<k\}$. Moreover, for every $k<\omega$
we obtain $b_k\equiv_c a_{n_k}$ and hence, $\models\varphi(c,b_k)$ for all $k<\omega$.
\eprova

\subsection{Weight and lowness} We have pointed out that supersimple theories of finite $\D$-rank are low. More generally, we will show that simple theories with bounded finite weight are also low. Recall the definition of pre-weight and weight:
\defi The pre-weight of a type $\tp(a/A)$, $\pwt(a/A)$, is the supremum of the set of all cardinals $\kappa$ for which there is an independent over $A$ sequence $(a_i:i<\kappa)$ such that $a\nind_A a_i$ for all $i<\kappa$. The weight of a type $\tp(a/A)$, denoted by $\w(a/A)$, is the supremum of the set of all pre-weights of the non-forking extensions of $\tp(a/A)$. \edefi

In a simple theory, by the local character of non-forking independence, every type has bounded (pre-)weight.

\defi We say that a simple theory has {\em bounded finite weight} if for any finite tuple of variables $x$ there is some natural number $n_{|x|}$ such that the weight of any type on $x$ (over any set of parameters) is bounded by $n_{|x|}$.\edefi

%We will say that the theory has a {\em bounded finite (pre-)weight} if there is a natural number $n<\omega$ such that every type has (pre-)weight less than $n$.

\rem Every simple theory of finite $\D$-rank has bounded finite weight as every complete type have bounded finite Lascar rank. Moreover, there are examples of simple non-supersimple theories all whose types on one variable have weight $1$ (e.g. $\mathrm{dp}$-minimal stable theories).  \erem

\prop Every simple theory of bounded finite weight is low. \eprop
\prova By Fact \ref{fact-type-def} it is enough to show that dividing is $\emptyset$-type-definable. Let $\varphi(x,y)\in L$ be a formula with $|x|=n$. By assumption there is some $k_n<\omega$ such that every complete type on $x$ has weight less than $k_n$. Firstly, we will check that for any tuple $a$ and any set $A$, $\varphi(x,a)$ divides over $A$ iff it divides over $A$ with respect to $k_n$. For this, consider a Morley sequence $(a_i:i<\omega)$ in $\tp(a/A)$; in particular, $\varphi(x,a_i)$ divides over $A$ for all $i<\omega$. Thus, for any $b$ such that $\models\varphi(b,a_i)$ with $i<\omega$ we have, $b\nind_A a_i$. As the $a_i$'s are $A$-independent, the set $\{\varphi(x,a_i):i<\omega\}$ must be $k_n$-inconsistent as otherwise we would obtain a type on $x$ over $A$ whose weight would be $\geq k_n$, a contradiction. Finally, it is clear that `$\varphi(x,y)$ divides over $z$ with respect to $k_{|x|}$' is $\emptyset$-type-definable on $y,z$. Hence, the result. \eprova

%\prova Assume there is a nonlow formula $\varphi(x,y)\in L$. Then for every $n<\omega$ there is a sequence $(a_i:i\leq n)$ exemplifying that $\varphi(x,y)$ divides $n+1$ times with respect to some $k_n$. Note that for every $n<\omega$, we may assume that $k_n$ is minimal and that $k_n<k_m$ for $n<m$, as otherwise a compactness argument yields the existence of a dividing chain of length $\omega_1$ contradicting simplicity. For all $n<\omega$, there is a Morley sequence $(a_n^l:l<\omega)$ in $\tp(a_n/a_j:j<n)$ witnessing that $\varphi(x,a_n)$ divides over $\{a_j:j<n\}$ with respect to $k_n$. But then for every $n<\omega$ and for any $b_n$ realizing the set $\{\varphi(x,a_n^l):l<k_n-1\}$, we have $b_n\nind_{(a_j:j<n)} a_n^l$ for all $l<k_n-1$, whence $\pwt(b_n/a_j:j<n)\geq k_n-1$, a contradiction.\eprova

\rem By \cite[Proposition 18]{cas-wag}, note that the proof above is also true for a simple theory with bounded finite weight for types on one variable. Moreover, the same is true if we require that for any cardinal $\lambda$, the weight of all types on one variable over any set of size $\lambda$ is bounded by a finite number.   \erem

%\rem In a simple theory of bounded finite weight, for every formula $\psi(x,y)\in L$ there is a $k<\omega$ such that a formula $\varphi(x,a)$ divides over $A$ iff it divides over $A$ with respect to $k$. \erem \prova Given a formula $\psi(x,a)$ which divides over $A$, let $k$ be the maximum $\w(b/A)$ among all $b$'s such that $\models\psi(b,a)$. It is clear that $\psi(x,a)$ divides over $A$ iff it divides over $A$ with respect to $k+1$.\eprova

By \cite{bue,shami}, as an immediate corollary we get the following result:

\cor In a simple theory of bounded finite weight, Lascar strong types and strong types coincide over any set of parameters. \ecor

%\subsection{Weight, $\omega$-categorical simple theories, and lowness}
To finish this section, we present a result on (pre-)weight in $\omega$-categorical simple theories. %In particular, we show that $\omega$-categorical simple theories of finite weight are low.

\lem\label{pre-weight} Let $T$ be an $\omega$-categorical simple theory. For a finite tuple $a$ and a finite set $A$, the type $\tp(a/A)$ has finite pre-weight. %That is, there is no infinite independent over $A$ sequence $(b_i:i\in I)$ such that $a\nind_A b_i$ for all $i\in I$.
\elem \prova Assume not, then for every $n<\omega$ there is an $A$-independent sequence $(b_i:i<n)$ with $a\nind_A b_i$ for all $i<n$. For a cardinal $\kappa$ big enough we consider a set of formulas $\Sigma(x_i:i<\kappa)$ expressing: \begin{center} $x_i\ind_A (x_j:j<i)$ and $a\nind_A x_i$ for all $i<\kappa$. \end{center} The first part is type-definable over $A$ by $\omega$-categoricity since $(x_i:i<\kappa)$ is $A$-independent iff every finite subsequence is. In addition, $a\nind_A x_i$ is $Aa$-definable again  by $\omega$-categoricity.
A compactness argument yields that $\Sigma(x_i:i<\kappa)$ is consistent and therefore, a realization of $\Sigma$ witnesses that $\tp(a/A)$ has pre-weight at least $\kappa$. As the argument works for any $\kappa$, this contradicts simplicity. \eprova

\rem By Fact \ref{Fact-indiscernibles} and Lemma \ref{pre-weight}, if a formula is nonlow in an $\omega$-categorical simple theory, then there is no $\emptyset$-independent sequence witnessing this. \erem

\section{A Lemma on the bounded closure operator}

This section is devoted to study the bounded closure operator in $\omega$-categorical theories. The results presented here are easy but we have not seen them in the literature.

The imaginary version of the next lemma was suggested to us by David M. Evans. However, we think convenient to present it in a hyperimaginary version. %This generalization to hyperimaginaries is not difficult but will require the following fact:

%\fet\label{Fact-hyp} Any $A$-bounded hyperimaginary is an equivalence class of a bounded type-definable over $A$ equivalence relation. More precisely, if $a_E$ is an $A$-bounded hyperimaginary and $p(x)=\tp(a)$, then $a_E=a_F$ where $F(x,y)$ is given by $\exists z(p(z)\wedge E(x,z)\wedge E(z,y))\vee x\equiv_A^\KP y.$\efet

\lem\label{Lemma-bdd} Assume $T$ is $\omega$-categorical, let $a$
be a finite tuple, and let $A$ be an arbitrary set (possibly of hyperimaginaries). Then, there is some
$e\in\mathfrak{C}^{eq}$ such that $$\bdd(e)=\bdd(a)\cap \bdd(A).$$\elem
\prova Let $h$ be a hyperimaginary such that $\dcl(h)=\bdd(a)\cap\bdd(A)$. By Neumman's Lemma we choose some $b\equiv_h a$ with $\bdd(a)\cap \bdd(b)=\bdd(h)$. Note that
$\bdd(a)\cap\bdd(b)=\bdd(a)\cap\bdd(A)$ and define the following relation
$$xyEuv \Leftrightarrow \bdd(x)\cap \bdd(y)=\bdd(u)\cap
\bdd(v).$$ It is obvious that $E$ is an $\emptyset$-invariant equivalence
relation and so, it is $\emptyset$-definable by $\omega$-categoricity. Let now
$e=(ab)_E$ and notice that $e\in\bdd(a)\cap \bdd(b)$. So it remains to check that
$\bdd(a)\cap \bdd(b)\sub \bdd(e)$. For this we consider the orbit of $h$ under $\Aut(\mathfrak C/e)$, denoted by $\mathcal O_e(h)$, and we check that $\mathcal O_e(h)$ is small, i.e., $|\mathcal O_e(h)|<|\mathfrak C|$. Since $\Aut(\mathfrak C/e)$ fixes $\bdd(a)\cap\bdd(b)$ setwise and $h\in\bdd(a)\cap\bdd(b)$, $\mathcal O_e(h)\sub \bdd(a)\cap\bdd(b)$. Then, each $e$-conjugate of $h$ is $a$-bounded and hence, each $e$-conjugate of $h$ is an equivalence class of a bounded type-definable over $a$ equivalence relation \cite[Proposition 15.27]{cas2}. Since there is just a bounded number of such equivalence relations, the orbit of $h$ under $\Aut(\mathfrak C/e)$ must be small.\eprova

Given a hyperimaginary $h$ we shall write $\acl^{eq}(h)$ to denote the set of imaginaries which are bounded over $h$, that is, $\acl^{eq}(h)=\bdd(h)\cap\mathfrak C^{eq}$. Recall that a hyperimaginary is said to be {\em quasi-finitary} if it is bounded over a finite tuple.

\cor\label{Cor-bdd1} Let $T$ be $\omega$-categorical. If a quasi-finitary hyperimaginary $h$ is bounded in some other hyperimaginary $h'$, then there is some finite set of parameters $A\sub\acl^{eq}(h')$ such that $h\in\bdd(A)$. \ecor \prova Assume $h\in\bdd(a)$ for some finite tuple $a$ and let $h'$ be a hyperimaginary such that $h\in\bdd(h')$. By Lemma \ref{Lemma-bdd} there is some imaginary  $e\in\mathfrak C^{eq}$ such that $\bdd(a)\cap\bdd(h')=\bdd(e)$. Hence, $e\in\acl^{eq}(h')$ and $h\in\bdd(e)$. \eprova

A simple theory {\em admits finite coding} if the canonical base of any finitary type is a quasi-finitary hyperimaginary  \cite[Chapter 6.1.3]{Wagner-book}. In particular, one-based theories and supersimple theories admit finite coding.

\lem\label{lem-bdd2} Let $T$ be an $\omega$-categorical simple theory. If $\cb(a/A)\subseteq\bdd(B)$ for some finite set $B$, then $\cb(a/A)$ is interbounded with an imaginary. \elem  \prova This is an immediate application of Lemma \ref{Lemma-bdd} since $\bdd(\cb(a/A))=\bdd(\cb(a/A))\cap\bdd(B)$. \eprova

\cor\label{Cor-bdd2}  An $\omega$-categorical simple theory which admits finite coding is supersimple. \ecor \prova Let $a$ be a finite tuple and let $A$ be an arbitrary set. By assumption and Lemma \ref{lem-bdd2}, there is some $e$ imaginary such that $\bdd(\cb(a/A))=\bdd(e)$, whence $e\in\acl^{eq}(A)$ and hence, $e\in\acl^{eq}(A_0)$ for some finite subset $A_0\sub A$. On the other hand, $a\ind_e A$ and so, $a\ind_{A_0} A$. \eprova

%\rem In fact the Corollary above shows that any type $\tp(a/A)$ which admits finite coding, does not fork over some finite subset $A_0\subseteq A$. \erem

\section{Main results}

In this section we shall prove the main results. We will investigate two approaches to the problem: via the stability of forking and via geometric properties of forking.

\subsection{Strong stable forking}  A simple theory has {\em strong stable forking} if whenever a type $\tp(a/B)$ forks over $A$, then there is a stable formula $\phi(x,y)\in L$ such that $\phi(x,b)\in\tp(a/B)$ forks over $A$. Observe that $A$ might not be a subset of $B$; if we add the requirement $A\subseteq B$ this corresponds to {\em stable forking}. Thus, strong stable forking implies stable forking.

\teo A countable $\omega$-categorical simple theory with strong stable forking is low. \eteo \prova Assume the ambient theory has strong stable forking, but suppose, towards a contradiction, that there is a nonlow formula $\varphi(x,y)\in L$. By Fact \ref{Fact-indiscernibles} there is a dividing chain $(a_i:i<\omega)$ witnessing that $\varphi(x,y)$ is nonshort and let $b$ be a realization of $\{\varphi(x,a_i):i<\omega\}$. Observe that $\tp(b/a_i)$ divides over $\{a_j:j<i\}$ for all $i<\omega$, so for each $i$ there is a stable formula $\psi_i(x,y)\in L$ such that $\psi_i(x,a_i)\in\tp(b/a_i)$ divides over $\{a_j:j<i\}$. By $\omega$-categoricity, there is just a finite number of formulas (up to equivalence) on $x,y$; thus, we may assume that all $\psi_i$'s are equivalent to some $\psi(x,y)\in L$. Hence, $\psi(x,y)$ is a stable formula which divides $\omega$ times and so it is not low, a contradiction. Hence, the result. \eprova

\rem In fact, for the proof above it is just necessary strong stable forking over finite sets. \erem

\subsection{The CM-trivial case}

Recall the definition of CM-triviality.

\defi A simple theory is CM-trivial if for all $a\in\mathfrak C^{eq}$, and for all sets of parameters $A\sub B$: if $\bdd(aA)\cap\bdd(B)=\bdd(A)$, then $\cb(a/A)\sub\bdd(\cb(a/B))$. \edefi

In the definition of CM-triviality we have to deal with the bounded closure operator since canonical bases are hyperimaginaries. By \cite[Corollay 3.5]{pal-wag}, in our context, each hyperimaginary is equivalent to a sequence of imaginaries and so, we may replace the bounded closure $\bdd$ in favour of the imaginary algebraic closure  $\acl^{eq}$.

%\begin{fact}\label{Fact-CM} A small simple CM-trivial theory eliminates hyperimaginaries. \end{fact}

\teo A countable $\omega$-categorical simple CM-trivial theory is
low.\eteo \prova We may assume that canonical bases are sequences of imaginaries and so, we may work in $T^{eq}$. Let $\acl$ denote the imaginary algebraic closure.

Suppose, towards a contradiction, that there is a nonlow formula $\varphi(x,y)\in L$. Then by Lemma
\ref{indisc-seq} there are some $c$ and some $c$-indiscernible
sequence $(a_i:i<\omega)$ such that for every $i<\omega$:
$\varphi(x,a_i)$ divides over $\{a_j:j<i\}$ and
$c\models\varphi(x,a_i)$. Now we prolong the sequence
to a $c$-indiscernible sequence $(a_i:i\leq\omega)$. Since
$\tp(a_\omega/a_i:i<\omega,c)$ is finitely satisfiable in
$\{a_i:i<\omega\}$ we have, $$a_\omega\ind_{(a_i:i<\omega)}c,$$ that is, $\cb(a_\omega/a_i:i<\omega,c)=\cb(a_\omega/a_i:i<\omega)$.

Let now $A=\acl(a_\omega c)\cap\acl(a_i:i<\omega,c)$. It follows that $c\in A=\acl(A)$ and that $\acl(A)=\acl(a_\omega A)\cap \acl(a_i:i<\omega,c)$. By CM-triviality we get $$\cb(a_\omega/A)\in\acl(\cb(a_\omega/a_i:i<\omega,c))=\acl(\cb(a_\omega/a_i:i<\omega)),$$ whence $\cb(a_\omega/A)\in\acl(a_i:i<\omega)$. Also, observe that $\cb(a_\omega/A)\in\acl(a_\omega c)$ and so, by Corollary \ref{Cor-bdd2} it is interalgebraic with a single imaginary element, say $e\in\mathfrak C^{eq}$. Hence, $a_\omega\ind_e A$; so, $a_\omega\ind_e c$. On the other hand, since $e$ is a single imaginary, there exists some $n<\omega$ such that $e\in\acl(a_i:i<n)$. But by $c$-indiscernibility observe that $\models\varphi(c,a_\omega)$ and that $\varphi(x,a_\omega)$ divides over $\acl(a_i:i<n)$, and so does over $e$; a contradiction. Hence, the result.\eprova

\preg The same proof will work without assuming CM-triviality if for all finite tuples $a,b$ and for every set $B$ with $b\in\bdd(B)$, there is some $\hat b\in\bdd(B)$ such that $b\in\bdd(\hat b)$ and $\cb(a/\hat b)\in\bdd(\cb(a/B))$, where $\hat b$ might be a quasi-finitary hyperimaginary. Is this true in general? \epreg

This question was already stated in \cite{pal-wag} where the author and Wagner observe that every theory satisfying this property would eliminate all hyperimaginaries if it eliminates finitary ones.

%Finally, we show that CM-trivial theories of finite weight have bounded finite weight.

%\prop Let $T$ be an $\omega$-categorical simple CM-trivial theory. Then $T$ have finite weight iff there is a finite bound on the pre-weight of types over finite sets. \eprop

\bibliographystyle{plain}

\end{document}